\documentclass[12pt,a4paper,reqno]{amsart}
\usepackage{amssymb}
\usepackage{dcpic,pictex}
\usepackage{multirow}
\usepackage{graphicx}
\usepackage{cite}

\usepackage{hyperref}
\newcommand\valor{0.49} 

\usepackage{hyperref}
\usepackage[main=english,french]{babel}

\newcommand{\Com}{\mathbb{C}}
\newcommand{\N}{\mathbb{N}}

\begin{document}


\title[Multidimensional scaling and visualization of patterns]{Multidimensional scaling and visualization of patterns in distribution of nontrivial zeros of the zeta-function}

\author{Jos\'{e} Tenreiro Machado$^1$}
\curraddr{$^1$Institute of Engineering, Polytechnic of Porto,  
     Department of Electrical Engineering, 
     Porto, Portugal}
\email{jtm@isep.ipp.pt}
\author{Yuri Luchko$^2$}
\curraddr{$^2$Beuth Technical University of Applied Sciences Berlin,  
     Department of  Mathematics, Physics, and Chemistry,  
     Luxemburger Str. 10,  
     13353 Berlin, Germany}
\email{luchko@beuth-hochschule.de}

\dedicatory{}
\keywords{Riemann's zeta-function, zeros of the zeta-function, distribution of zeros, complex systems, Lorentzian metric, multidimensional scaling  algorithm, periodical patterns}

\begin{abstract}
In this paper, we analyze  the nontrivial zeros of the Riemann zeta-function using the multidimensional scaling (MDS) algorithm and computational visualization features. The nontrivial zeros of the Riemann zeta-function as well as the vectors with several neighboring zeros are interpreted as the basic elements (points or objects) of a data set. Then we employ a variety of different  metrics, such as the Euclidean and Lorentzian ones, to calculate the distances between the objects. The set of the calculated distances is then processed by the MDS algorithm that produces the loci, organized according to the objects features. Then they are analyzed from the perspective of the emerging patterns. Surprisingly, in the case of the Lorentzian metric, this procedure leads to the very clear periodical structures both in the case of the objects in form of the single nontrivial zeros of the Riemann zeta-function and in the case of the vectors with a given number of neighboring zeros. The other tested metrics do not produce such periodical structures, but rather chaotic ones. In this paper, we restrict ourselves to  numerical experiments and the visualization of the produced results. An analytical explanation of the obtained periodical structures is an open problem worth for investigation by the experts in the analytical number theory. 
\end{abstract}

\maketitle

\section{Introduction}\label{sec:1}

In his seminal book \cite{Euler:1748}, Leonard Euler not just settled the foundations of the modern calculus, he also provided a series of inspiring ideas that led to several new theories and important results in different branches of mathematics. In particular,  Euler introduced a series that corresponds to the Riemann zeta-function of the real argument
\begin{equation}
\label{zeta-real}
\zeta(x) := \sum_{n=1}^{+\infty} \frac{1}{n^x}, \ x>1
\end{equation}
and derived the explicit formulas for some of its values including $\zeta(2n),\ n\in \N$. In the same book, Euler presented one of the most beautiful and important formulas of the number theory that connects the zeta-function with the prime numbers:
\begin{equation}
\label{Euler}
\prod_{k=1}^{+\infty} \left(1-\frac{1}{p_k^x}\right)^{-1} = \sum_{n=1}^{+\infty} \frac{1}{n^x},
\end{equation}
where $p_k=2,3,5,\dots$ are the prime numbers ordered by their values. 

Several mathematicians and primarily the great Russian mathematician P.L. Tche\-bi\-chef in \cite{Tchebichef:1852} employed the Euler formula \eqref{Euler} restricted to the real numbers $x>1$ for investigation of the prime numbers distribution. However, the crunch step in the right direction was done by Bernhard Riemann in \cite{Riemann:1859}, where he first extended the convergence domain of \eqref{zeta-real} to a part of the complex plane 
\begin{equation}
\label{zeta-com}
\zeta(z) := \sum_{n=1}^{+\infty} \frac{1}{n^z}, \ \Re(z)>1
\end{equation}
and then introduced his famous zeta-function as an analytical continuation of the series \eqref{zeta-com} into the domain $\Com \setminus \{ 1\}$. The aim of Riemann was to prove the hypothesis of Gauss-Legendre that says that the number $\pi(x)$ of the prime numbers that do not exceed $x$ satisfies the following asymptotic relation:
\begin{equation}
\label{Gauss}
\pi(x) \sim \mbox{Li}(x)  \ \mbox{as} \ x\to +\infty,\ \mbox{Li}(x):= \int_2^x \frac{dt}{\ln(t)} \sim \frac{x}{\ln(x)}  \ \mbox{as} \ x\to +\infty,
\end{equation}
where $\mbox{Li}$ stands for the integral logarithm. In \cite{Riemann:1859}, Riemann established a connection between the formula  \eqref{Gauss} and the distribution of the non-trivial zeros of his zeta-function (i.e., the zeros not equal to $-2,\, -4,\, -6,\dots$). Moreover, he derived an explicit formula for the function $\pi(x)$ in terms of the $\mbox{Li}$-function and the non-trivial zeros of the zeta-function (in slightly different form):
\begin{equation}
\label{Riem_1}
\pi(x) = \mbox{Li}(x) - \frac{1}{2} \mbox{Li}(x^{1/2})+\sum_{\rho_i} \mbox{Li}(x^{\rho_i}) + \mbox{some small terms},
\end{equation}
where $\rho_i,\ i=1,2\dots$ are nontrivial zeros of the zeta-function. 
The hypothesis of Gauss-Legendre follows from the Riemann formula \eqref{Riem_1} when the condition $\Re(\rho_i)<1$ is fulfilled for all nontrivial zeros of the zeta-function.  Riemann by himself could prove in \cite{Riemann:1859} that all  nontrivial zeros of the zeta-function lay in the so-called critical strip: $0\le \Re(\rho_i) \le 1$. The sharp inequality $\Re(\rho_i)<1$ and thus  the Gauss-Legendre hypothesis \eqref{Gauss} was proved by Hadamard in \cite{Hadamard:1896} and by de la Vall$\acute{\mbox{e}}$e Poussin in \cite{Poussin:1896} independently each from other  about 40 years after publication of  the Riemann's famous paper \cite{Riemann:1859}.

Another and in fact much more important contribution of Riemann to the problem of prime numbers distribution was formulation of his famous hypothesis (known as the Riemann hypothesis, or abbreviated as RH) that says that all nontrivial zeros of the zeta-function lay on the straight line with the real part one half: $\Re(\rho_i) = \frac{1}{2}$. The RH is not proved until now despite of numerous efforts of several generations of mathematicians. It is one of the so called millennium problems established by the Clay Mathematics Institute in 2000 (\url{http://www.claymath.org/millennium-problems}). 

The impact of the RH to the number theory cannot be overestimated. As shown already by Riemann, its validity would lead to the following asymptotic formula for the distribution of the prime numbers:
\begin{equation}
\label{Riem_3}
\pi(x) = \mbox{Li}(x) + O(x^{1/2}\log(x)),\ x\to +\infty.
\end{equation}
In fact, the statement  \eqref{Riem_3} is equivalent to the RH. Moreover, the RH can be reformulated not only in terms of asymptotic behavior of the prime numbers distribution, but also as some important statements in different fields of science as, e.g., theory of the finite groups transformations, probability theory, functional analysis, and even quantum mechanics. 


%



The RH is probably the most important and famous unsolved mathematical problem. Nowadays there are so many indicators that the RH is true that everybody believes in its validity. For example, Hardy proved in \cite{Hardy:1914} that infinitely many nontrivial zeros of the zeta-function lay on the critical line $\Re(z) = 1/2$. Levinson improved this result in \cite{Levinson:1974} and showed that more than one third of the nontrivial zeros lay on the critical line. In \cite{Lune:1986}, van de Lune, te Riele, and  Winter verified that the first one and half milliards of the nontrivial zeros of the zeta function lay on the critical line.  For further details regarding the RH we refer the interested readers to \cite{Hardy:1921}, \cite{Selberg:1942}, \cite{Titchmarsh:1987} for discussion of the classical results and to \cite{Bombieri:2006} and \cite{Karatsuba:1992} for surveys of the recent results.  

Along with the analytical results mentioned above, the numerical methods for calculation of the zeta-function and its nontrivial zeros with an arbitrary precision were worked out (see \cite{Odlyzko:1988}, \cite{Moree:2018} and references therein). Nowadays, there exist numerous  tables of numerical values of millions of the nontrivial zeros of the zeta-function calculated with a high precision (see, e.g., the website (\url{https://www.lmfdb.org/zeros/zeta/}). Moreover, the CAS Mathematica provides a special command \texttt{ZetaZero}[k] that calculates the $k$th nontrivial zero of the zeta-function with an arbitrary precision. 

It is 	widely recognized that the nontrivial zeros of the zeta-function are distributed rather irregular and random. All the more surprising is the fact that we could found out a certain periodicity in the distribution of the Lorentzian distances between the nontrivial zeros of the zeta-function as well as between the vectors containing several neighboring  nontrivial zeros of the zeta-function. 

Stemming from these ideas the manuscript is organized as follows. Section \ref{sec:2} introduces the fundamental concepts regarding the MDS algorithm. The role of distance between objects and the properties of six metrics, namely the  angular Arccosine, Jaccard, Chebyshev, Euclidean, Canberra and Lorentzian, are also discussed. Section \ref{sec:3} formulates the problem and discusses the results. The data-set and the construction of vectors with consecutive values of zeros of the zeta function is first described. Then the performance of the MDS for six test distances and several distinct approaches is assessed. The results with the  Lorentzian reveal a periodicity and their properties are further analyzed in numerical terms. Finally, Section \ref{sec:4} presents the main conclusions.

\section{Multidimensional scaling algorithm}\label{sec:2}

MDS is an iterative algorithmic technique \cite{Torgerson:58,Kruskal:78} that visualizes in a $n$-dimensional space objects initially described in a $m$-dimensional space, where in general $m\geq n$. The main idea is to achieve their visualization by means of a graphical representation  where the objects are represented by points. For that purpose MDS uses the concept of distance for comparing the objects.

A function $d$ gives a distance between two objects $\xi_1$ and $\xi_2$ if satisfies the three axioms \cite{Deza:09}:

\begin{subequations} \label{eq:metrics}
\begin{align}	
  		d(\xi_1,\xi_2)&= 0, \ \text{if} \ \xi_1 = \xi_2, \text{identity},\\  	
  		d(\xi_1,\xi_2)&= d(\xi_2,\xi_1), \text{symmetry},\\  
  		d(\xi_1,\xi_2)&\leq d(\xi_1,\xi_3) + d(\xi_2,\xi_3), \text{triangle inequality}.
\end{align}
\end{subequations}

These axioms admit using different functions each with its own pros and cons.  Therefore, for a given set of objects, in general it is a good strategy to test in advance a number of distances and, based on the results, to select those that reflect more adequately the phenomenon under analysis \cite{Machado:2021}. 

Let us denote by $\xi_{i}=\left[x_{i}(1),\ldots,x_{i}(m)\right]$ and $\tilde{\xi}_{i}=\left[\tilde{x}_{i}(1),\ldots,\tilde{x}_{i}(n)\right]$ the $i$-th object in the original and approximation spaces, respectively, and consider a given distance $d$. The process is initialized by calculating a $N\times N$ symmetrical matrix $\left[d(\xi_i,\xi_j)\right]$ of  distances between the objects $\xi_i$ and $\xi_j$. The MDS algorithm produces the new set of objects $\tilde{\xi_i}$, $i=1,\ldots,N$,  that minimizes a given optimization index, often called stress $S$. The problem is, therefore, converted to the numerical minimization of $S$ and commonly users adopt 2 or 3 dimensions, because this allow a direct portraying.  When opting for 3 dimensions we have a slightly superior approximation than just for 2, but that comes with some extra cost, namely with the need for rotation, shift and amplification for obtaining a good visualization. Obviously, we can analyze each set coordinates separately \cite{Machado:2020c}, but in most practical cases that is not necessary.

The quality of the approximation through the objects $\tilde{\xi}$ can be assessed by means of the Sheppard and stress plots \cite{Sammon:69,Borg:05}. The Sheppard diagram draws $d_{ij}$ versus $\tilde{d}_{ij}$ and, consequently, a low/high scatter means a good/poor match. The stress diagram plots $S$ versus $n$, which gives a monotonic decreasing curve. We verify with a significant reduction of $S$ at the initial values and $n=3$ is in general a good compromise between accuracy and feasibility of portraying. The process is concluded by  the user and consists of analyzing  the MDS locus.  We must note that the axes have no physical meaning and there is no interpretation of good/bad for high/low values of $\tilde{\xi}$. Instead, users must read the MDS results in the perspective of the clusters and patterns formed by the objects $\tilde{\xi}$ since they reflect some relationship embedded in the original objects under the light of the distance $d$. As mentioned before, it is advisable to test several distances and to interpret the corresponding MDS locus for selecting the `best'.

Hereafter, we consider a test-bed of six distances, namely the angular Arccosine, Jaccard, Tchebichef, Euclidean, Canberra, Lorentzian, given by:

\begin{subequations} \label{eq:distances}
\begin{align}
d^{Ac}(\xi_i,\xi_j) & =\arccos \displaystyle  \frac{{\displaystyle  \sum_{k=1}^{m}x_{i}\left(k\right)x_{j}\left(k\right)}}{\sqrt{\displaystyle  \sum_{k=1}^{m}x_{i}\left(k\right)^{2}\sum_{k=1}^{m}x_{j}\left(k\right)^{2}}},\\
d^{Ja}(\xi_i,\xi_j)  & ={\displaystyle \frac{\displaystyle  \sum_{k=1}^{m}x_{i}(k)x_{j}(k)}{\displaystyle  \sum_{k=1}^{m}x_{i}(k)^{2}+\sum_{k=1}^{m}x_{j}(k)^{2}-\sum_{k=1}^{m}x_{i}(k)x_{j}(k)}},\\
d^{Tc}(\xi_i,\xi_j)  & =\underset{k}{\max}\left(\left|x_{i}(k)\right|-x_{j}(k)\right),\\
d^{Eu}(\xi_i,\xi_j)  & ={\displaystyle \sqrt{\sum_{k=1}^{m}\left[x_{i}\left(k\right)-x_{j}\left(k\right)\right]^{2}}},\\
d^{Ca}(\xi_i,\xi_j)  & ={\displaystyle \sum_{k=1}^{m}\frac{\left|x_{i}\left(k\right)-x_{j}\left(k\right)\right|}{\left|x_{i}\left(k\right)\right|+\left|x_{j}\left(k\right)\right|}},\\
d^{Lo}(\xi_i,\xi_j)  & ={\displaystyle \sum_{k=1}^{m}\ln\left[1+\left|x_{i}\left(k\right)-x_{j}\left(k\right)\right|\right]},
\end{align}
\end{subequations}

\noindent
where $x_{i}\left(k\right)$ and $x_{j}\left(k\right)$, $i,j=1,\ldots,N$, denote the $k$-th components of the $i$-th and $j$-th objects, respectively. 

The Arccosine distance is not sensitive to amplitude and just provides a measure of the angle between two vectors. The Jaccard distance measures the dissimilarity between two sample sets and is useful for comparing observations with categorical variables.  The  Tchebichef distances and Euclidean are special cases of the Minkowski distance $d^{Mi}=\left[\sum_{k=1}^{m}\left|x_{i}\left(k\right)-x_{j}\left(k\right)\right|^{p}\right]^{1/p}$ for $p\rightarrow\infty$ and $p=2$, respectively. The Canberra distance is a weighted version of the Manhattan that we obtain when substituting $\left|x_{i}\left(k\right)-x_{j}\left(k\right)\right|$ by $\frac{\left|x_{i}\left(k\right)-x_{j}\left(k\right)\right|}{\left|x_{i}\left(k\right)\right|+\left|x_{j}\left(k\right)\right|}$. Similarly, the Lorentzian distance adjusts the comparison of small and large values by means of the $\log(\cdot)$ function. Therefore, the Canberra and Lorentzian distances adapt well in cases involving small and large numerical values. Finally, we calculate the MDS using the Matlab classical multidimensional scaling command \texttt{cmdscale} \cite{Martinez:05}.

\section{Periodic patters and their visualization}\label{sec:3}

The list of consecutive zeros $x(k)$, $k=1,\ldots$, of the zeta function was retrieved from \url{https://www.lmfdb.org/zeros/zeta/} . Therefore, the $i$-th object in the MDS process consists of the $i$-th vector $\xi_{i}=\left[x_{i}(1),\ldots,x_{i}(m)\right]$. Moreover,  for constructing two consecutive vectors with values selected from the list, two approaches are designed, namely $\mathcal{A}_{1}$: with disjoint values and $\mathcal{A}_{2}$: sharing common values excepting the first. More specifically, for two consecutive objects, $i$-th  and $i+1$-th, we have $\mathcal{A}_{1}$:  $\xi_{i+1}=\left[x_{i}(m+1),\ldots,x_{i}(2m)\right]$ and $\mathcal{A}_{2}$:  $\xi_{i+1}=\left[x_{i}(2),\ldots,x_{i}(m+1)\right]$.

The MDS 3 - dimensional loci in Fig. 1 are constructed from the list with the first $10,000$ zeros using the six distances (\ref{eq:distances}) and $m=10$ for approach $\mathcal{A}_{1}$, yielding  $N=1,000$. The vertical color bar indicates the order of the vector $i=1,\ldots, N$.


\begin{figure}[h]\label{fig:MDS1}
\begin{center}
\includegraphics[width=\valor\linewidth]{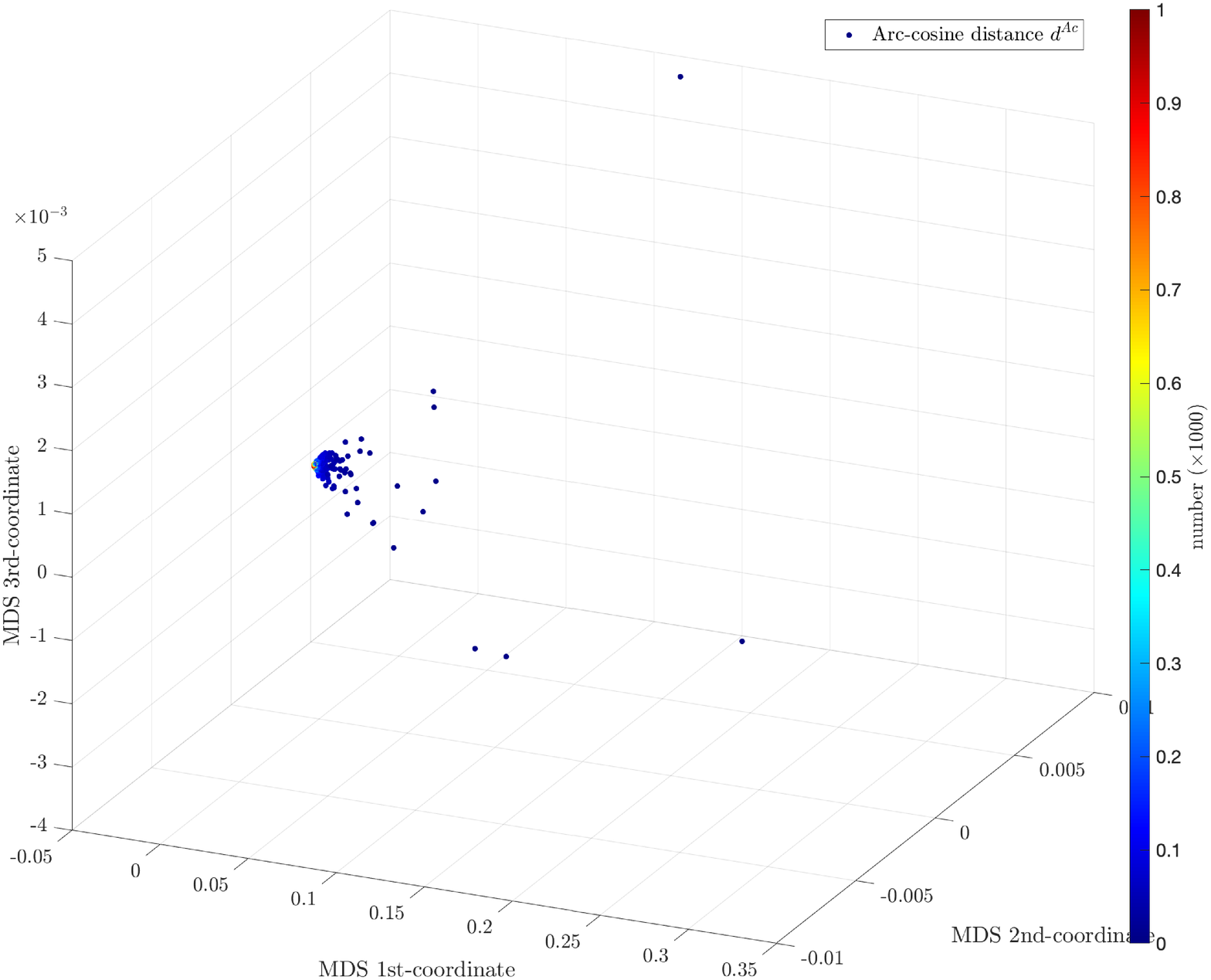}
\includegraphics[width=\valor\linewidth]{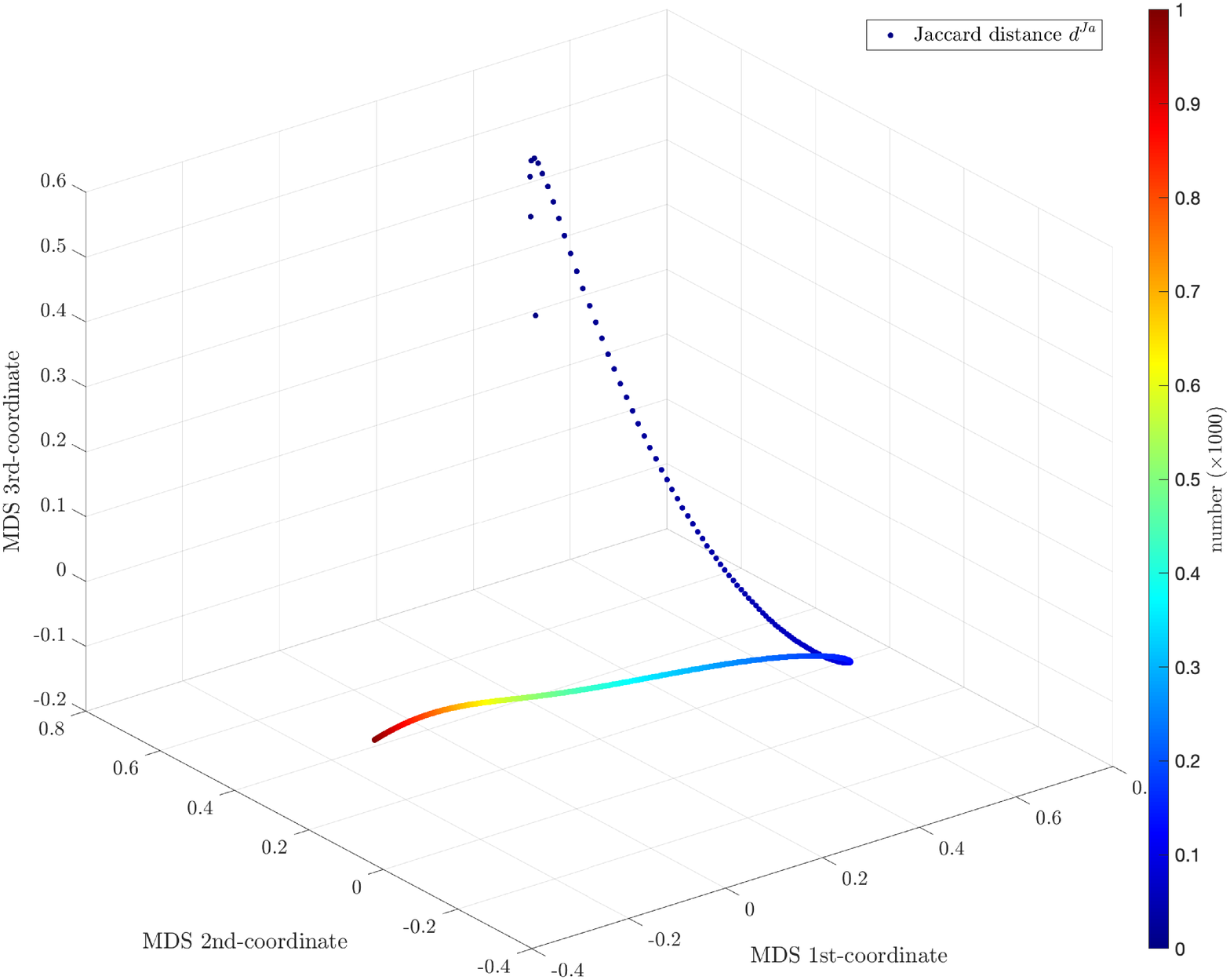}\\
\includegraphics[width=\valor\linewidth]{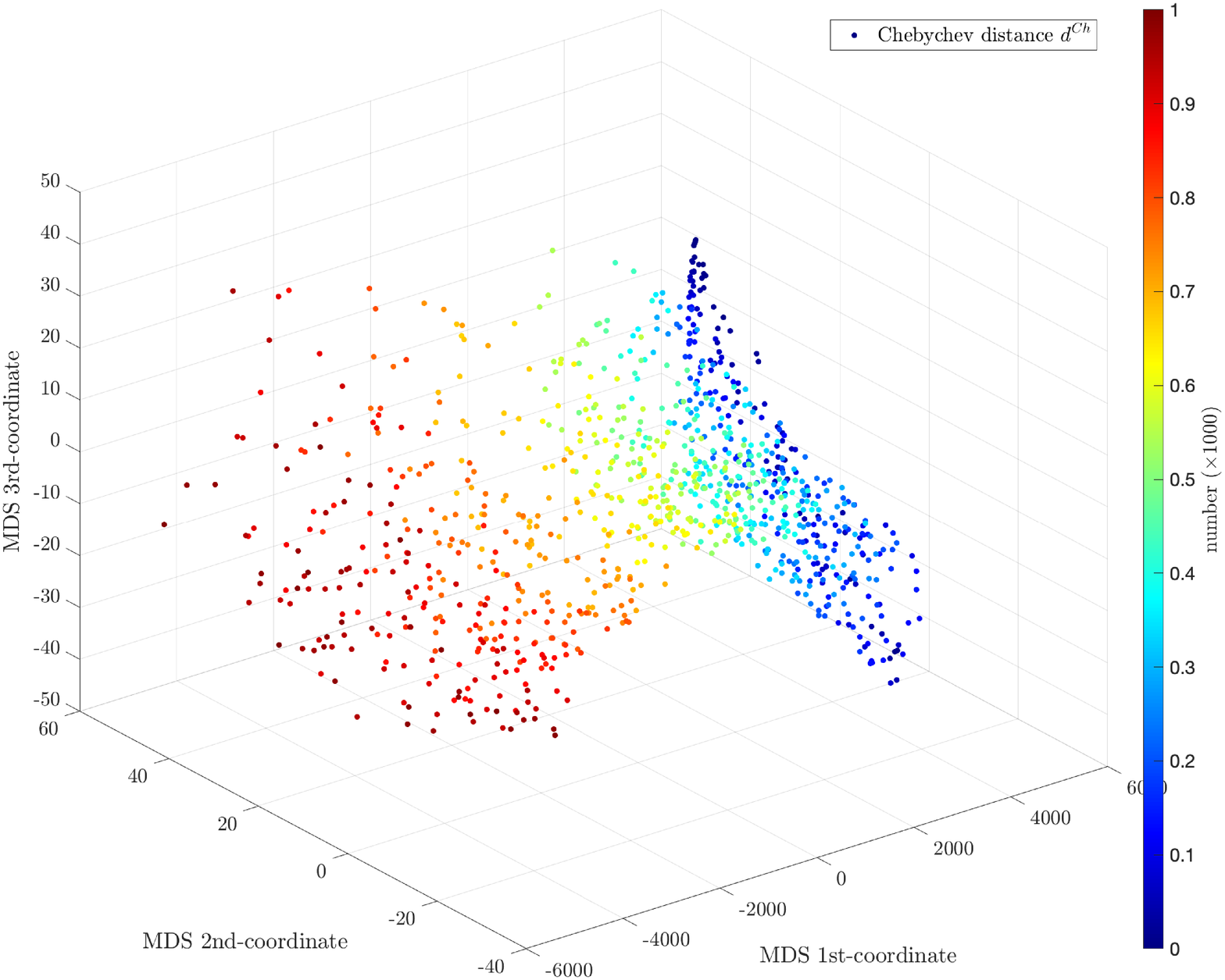}
\includegraphics[width=\valor\linewidth]{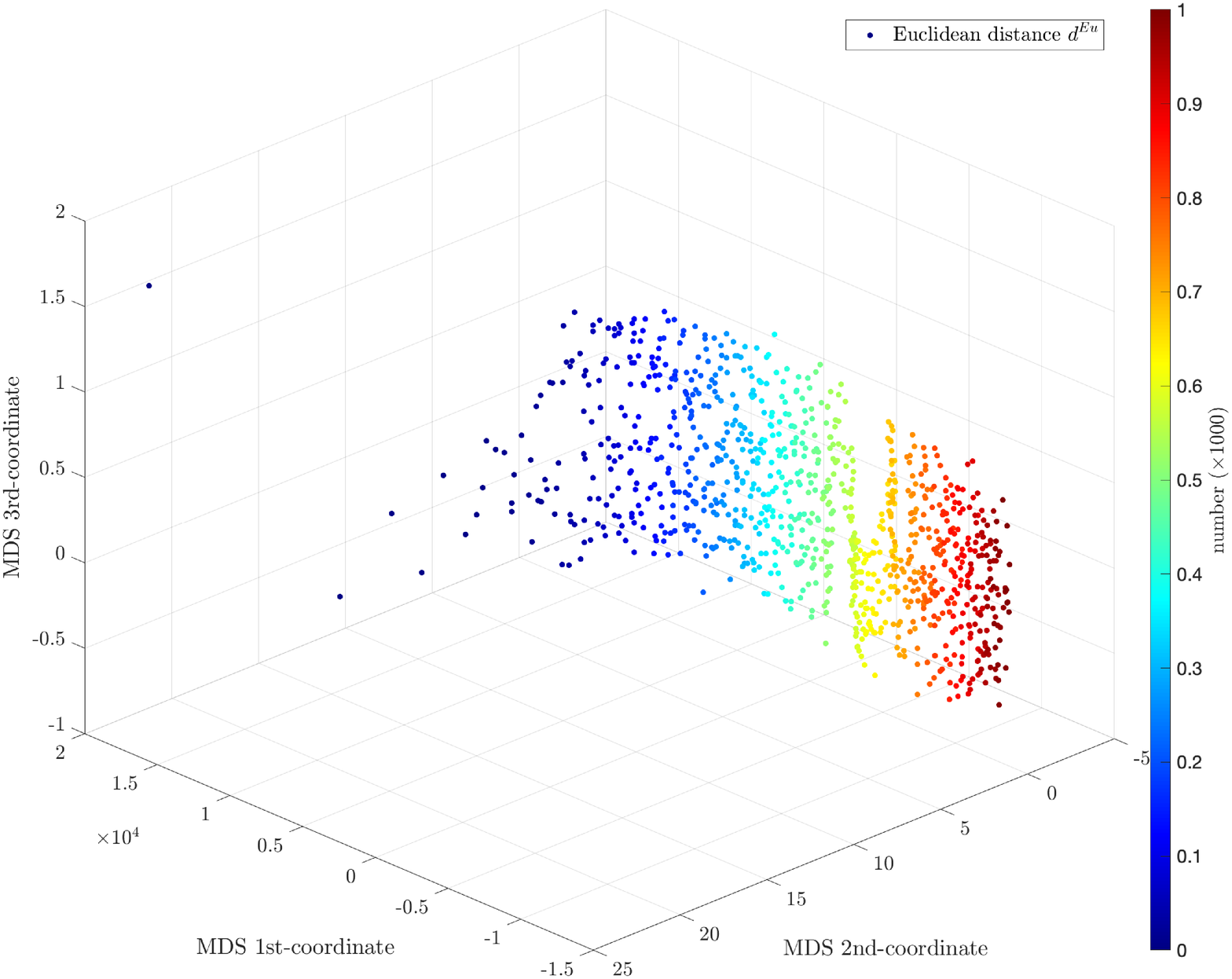}\\
\includegraphics[width=\valor\linewidth]{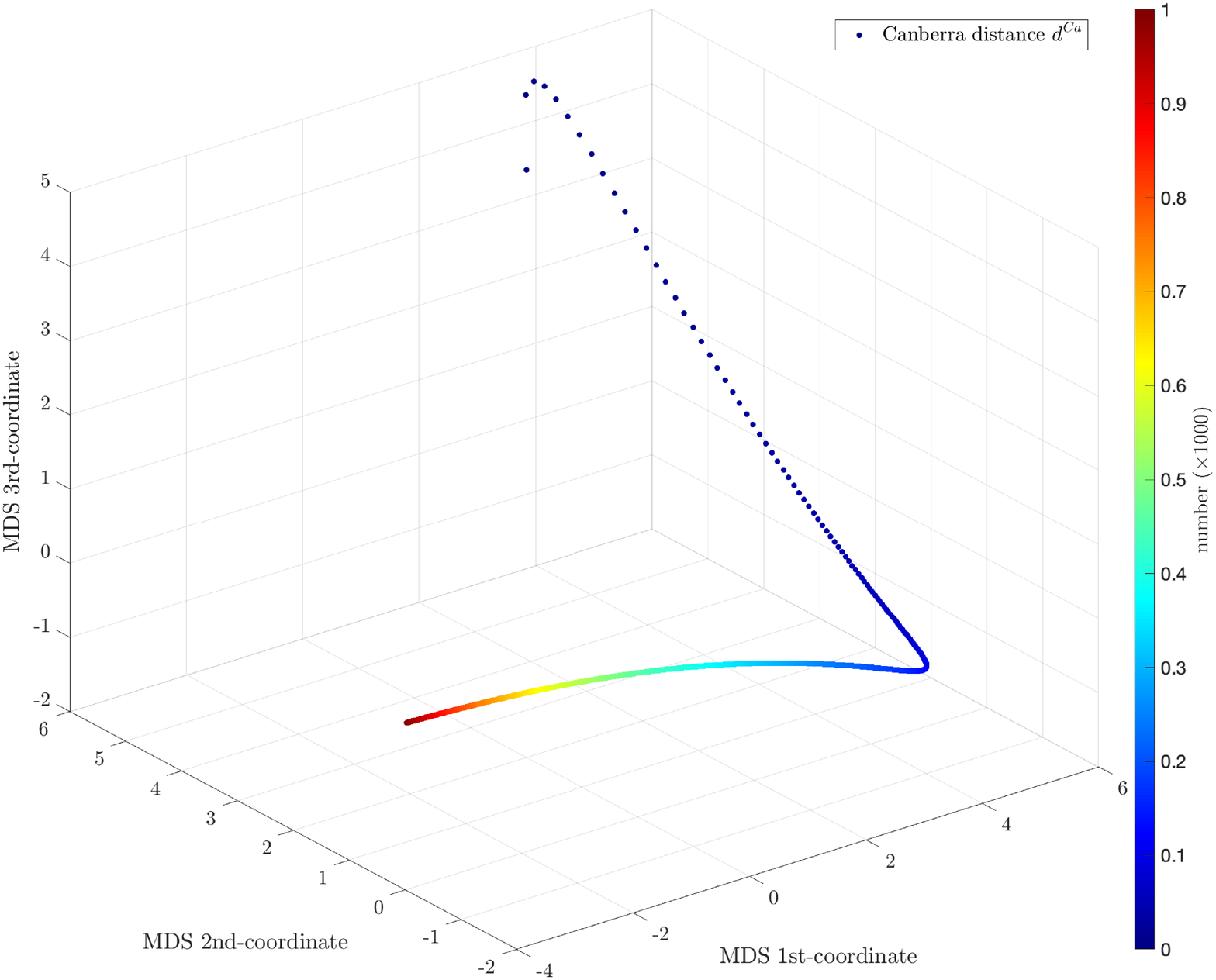}
\includegraphics[width=\valor\linewidth]{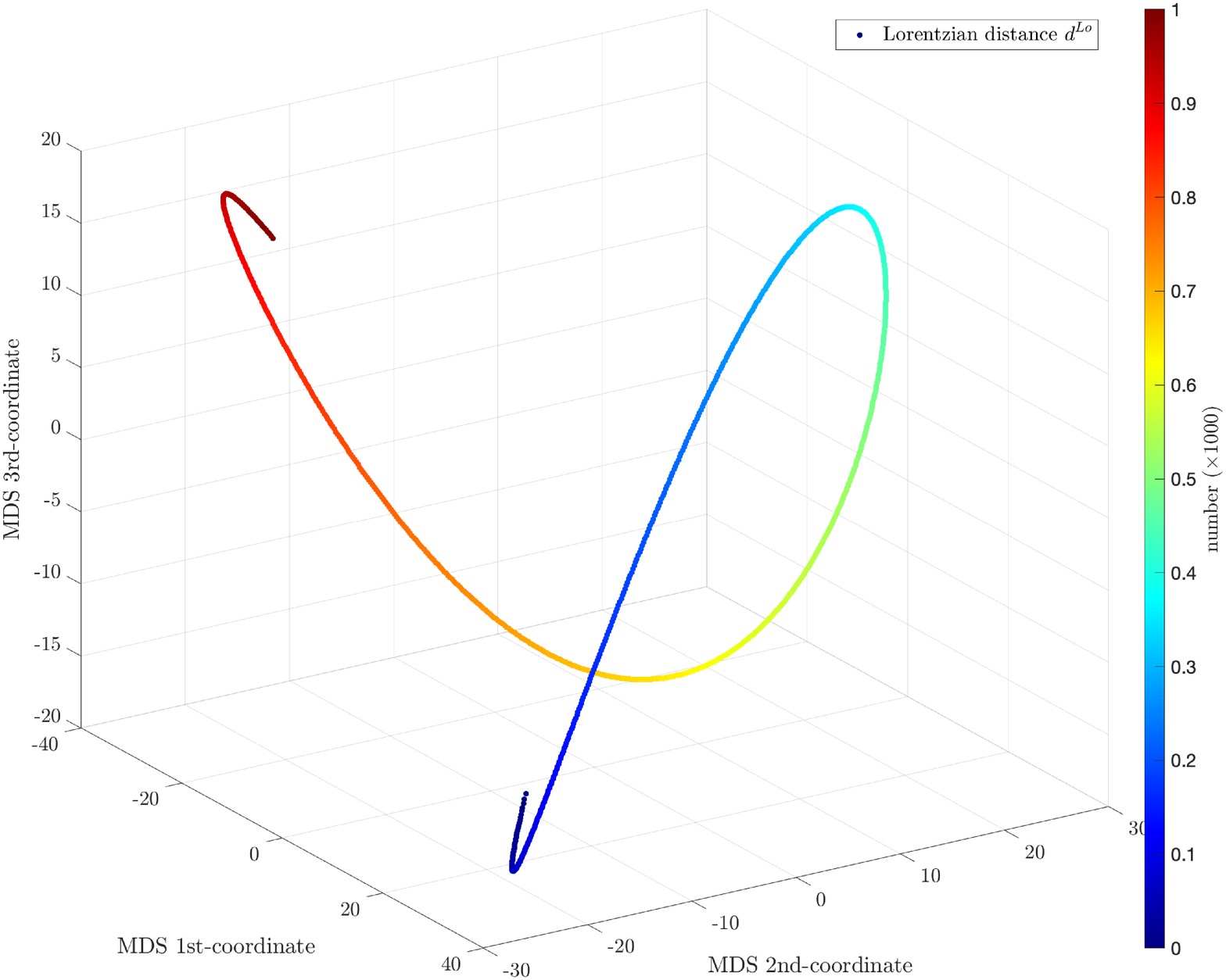}
\caption{The 3-dim MDS loci for $N=1,000$ vectors using the six distances (\ref{eq:distances}), for $m=10$ and approach $\mathcal{A}_{1}$.}
\end{center}
\end{figure}

We verify that $d^{Ac}$ is not capable of revealing any relationship. The pairs of distances $d^{Ja}$ and $d^{Ca}$, on one hand, and $d^{Tc}$ and $d^{Eu}$, on the other hand, show some kind of relationship, but without a clear pattern. However, $d^{Lo}$ shows clearly a periodic behavior. Exploring the approach $\mathcal{A}_{2}$ yields the same type of results and therefore the MDS plots are not included for being parsimonious. Moreover, repeating the analysis for $d^{Lo}$ and other values of $m$ gives also the MDS plots with a periodic-like behavior.

Figure 2 shows the MDS components, that is, the values of $\tilde{\xi}_{i}(p) $ versus $i$, $p=1,\ldots,10$,  for the distance $d^{Lo}$, when $m=10$ under the approach $\mathcal{A}_{1}$. We verify clearly an evolution close to a sinusoidal as anticipated from the 3-dim MDS locus. A numerical approximation with $\tilde{\xi}_{i}(p) \thickapprox A_p \sin\left(\omega_p i+\phi_p\right)$, $A_p,\omega_p\in \mathbb{R}^+$ and $ \phi_p\in \mathbb{R}$, yields the plots of Fig. 3. We observe a clear variation with $p$, with particular emphasizes in the power law  and  linear relations $ A_p \thicksim p ^{-0.5}$ and $ \omega_p \thicksim p$. 

In a critical analysis we can say that these results are merely for (i) vectors based on the zeta function zeros, (ii) numerical, and (iii) using MDS with the Matlab classical multidimensional scaling command \texttt{cmdscale}. Therefore, a more complete analytic analysis is needed to give a more solid basis. Nonetheless,  these results are relevant to explore new directions of research regarding the RH.

\begin{figure}[h]\label{fig:MDS2}
\begin{center}
\includegraphics[width=\linewidth]{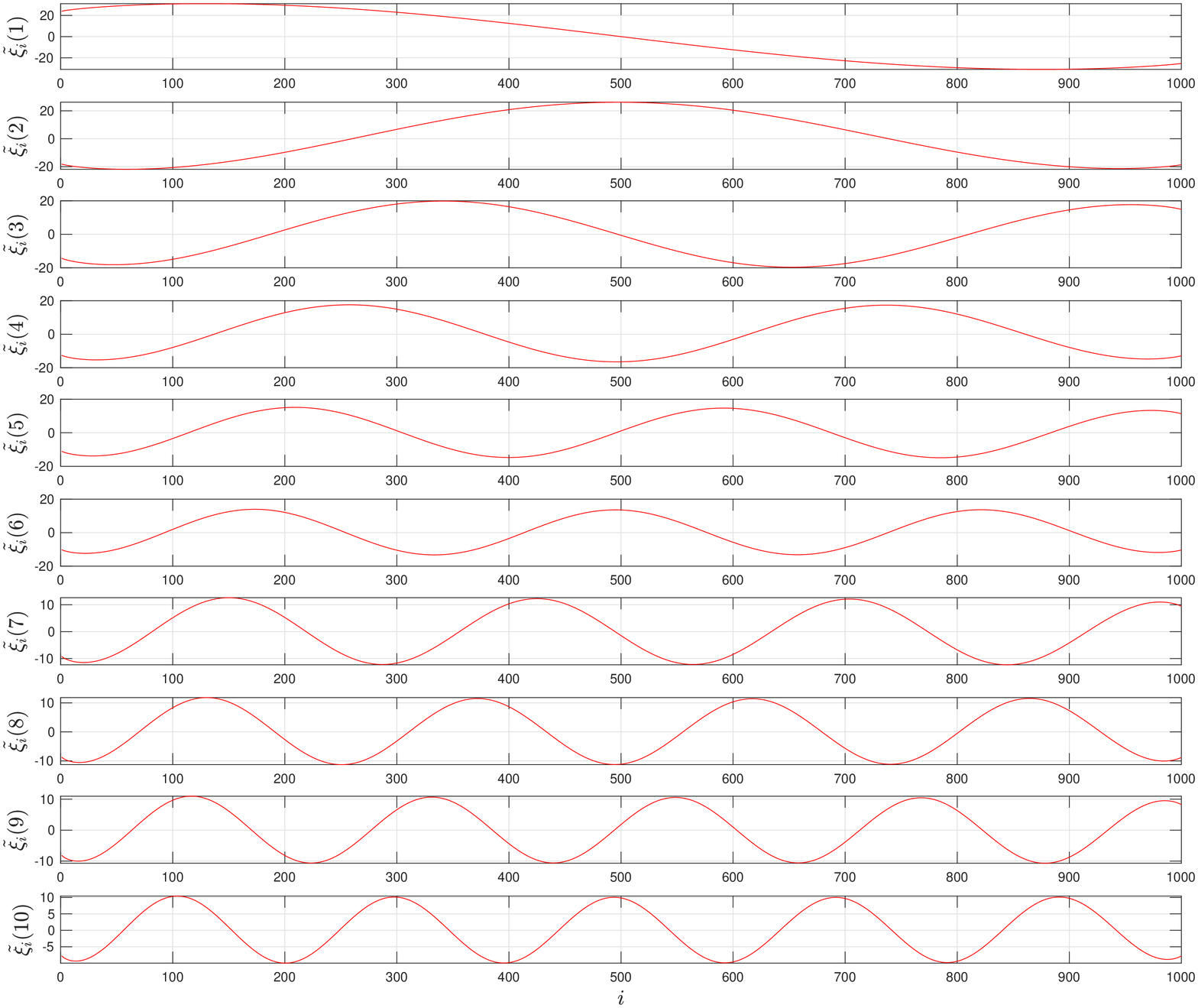}
\caption{Evolution of $\tilde{\xi}_{i}(p) $ versus $i$, $p=1,\ldots,10$, with the distance $d^{Lo}$, $m=10$ and $\mathcal{A}_{1}$.}
\end{center}
\end{figure}

\begin{figure}[h]\label{fig:MDS3}
\begin{center}
\includegraphics[width=\linewidth]{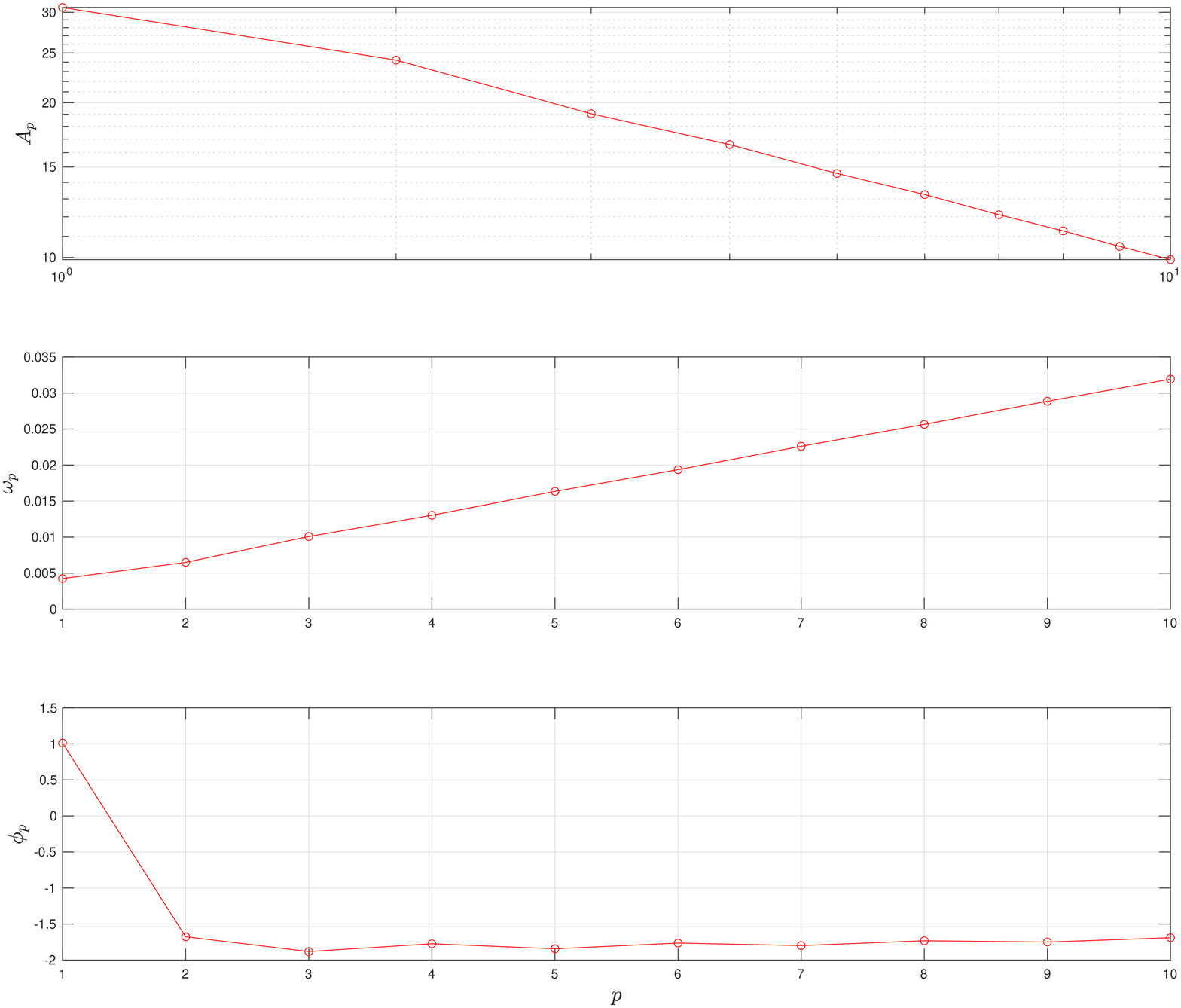}
\caption{Variation of $A_p$, $\omega_p$ and $\phi_p$ versus $p$, for the approximation $\tilde{\xi}_{i}(p) \thickapprox A_p \sin\left(\omega_p p+\phi_p\right)$, with the distance $d^{Lo}$, $m=10$ and $\mathcal{A}_{1}$.}
\end{center}
\end{figure}

\section{Conclusions}\label{sec:4}
 
This paper explored a possible perspective of the RH using the MDS computational technique. The MDS allows the visualizing possible relationships embedded into a data-set, which in our case consisted of the non-trivial zeros of the zeta function. The MDS requires the adoption of an appropriate distance, capable of capturing the characteristics of the phenomena to be unraveled. For that purpose six metrics were tested in conjunction with the classical MDS, yielding a very clear result for the Lorentzian distance. Indeed, a periodic behavior was found in the MDS representation of the vectors of consecutive zeros. A consequent numerical analysis of the periodicity yields also a power law and a linear evolution of the parameters of the sinusoidals exhibited in the MDS coordinates. These results represent  simply based on a computational and algorithmic approach, but may be a good starting point for new advances regarding the long standing RH problem.

\bibliographystyle{unsrt}
\bibliography{ref_V3}







\end{document}